# Effects of Axial Load on the Location of a Combined Null Point in Energy Piles

Arash Saeidi Rashk Olia[1], Dunja Perić[2]

[1] Department of Civil Engineering, Kansas State University, 1701C Platt St., Manhattan, KS 66506-5000; E-mail: saeidi@ksu.edu
[2] Department of Civil Engineering, Kansas State University, 1701C Platt St., Manhattan, KS, 66506-5000; E-mail: peric@ksu.edu

**ABSTRACT**

Soil structure interaction in energy piles has not yet been understood comprehensively. One of the important underlying issues is the location of a zero-displacement point, known as the null point. This study investigates how the location of a combined null point in fully floating energy piles is affected by relative magnitudes of thermal and mechanical loads. Analytical solutions are used to address four different loading scenarios including compressive and tensile loads, and heating and cooling. It was found that the location of the combined null point coincides with the location of the maximum magnitude of the axial stress induced by thermo-mechanical load. Furthermore, while a thermal null point is always present in fully floating energy piles, the combined null point is absent in the case of a small magnitude of thermal load as compared to the mechanical load. With increase in the relative magnitude of thermal load as compared to the mechanical load the combined null point emerges at different locations for different load combinations. In all cases it subsequently moves towards the mid-length of the pile with increase in the relative magnitude of thermal load.

**INTRODUCTION**

Energy piles are multifunctional deep foundations. They transfer superstructure loads to the ground while simultaneously enabling exchange of thermal energy between the superstructure and shallow subsurface. Consequently, energy piles are sustainable foundations that use a low enthalpy geothermal energy for a supplemental space heating and cooling, thus improving the energy efficiency of buildings while decreasing the emission of carbon dioxide into the atmosphere. The great potential for environmental and economic benefits stems from utilization of geothermal energy as a source of renewable energy, thus resulting in a worldwide popularity of geothermal piles (Ghasemi-fare and Basu 2016).

Due to restraints imposed by the surrounding soil, temperature change of an energy pile induces displacement, strain and stress in the pile (Bourne-Webb et al.2009, Amatya et al. 2012, Bourne-Webb et al. 2013, Perić et al. 2020). Tendency of a heated energy pile to expand and a cooled pile to contract is the main cause of induced thermal stresses. Furthermore, the soil that surrounds piles induces compressive stress in a heated energy pile and tensile stress in a cooled energy pile. Consequently, prediction of thermally and mechanically induced displacements, strains and stresses is vitally important for design of energy piles. Finite element models and analytical solutions showed that the combined null point emerges in semi floating energy pile tip, while the combined null point in perfectly end bearing pile is always at the pile tip, regardless of the load combination. (Saeidi Rashk Olia and Perić 2021 a & b)

This study focuses on evaluation of effects of relative magnitude of mechanical load with respect to thermal load on the location of a combined null point, which is a location of zero displacement. It is because semi floating and fully floating piles expand at both ends when heated and vice versa when cooled that there is no displacement at all at a certain location that is known as a null point. For example, in nearly fully floating piles the thermal null point is located at the mid-length of the pile.

**ANALYTICAL MODEL**

Analytical solutions for thermo-mechanical soil structure interaction in a single energy pile were presented by Cossel (2019), Iodice et al. (2020), and Perić, et. al. (2020). They all assumed that the soil pile interface remains in elastic regime while pile obeys a thermo-elastic constitutive law. These assumptions were confirmed by Laloui et al. (2006), Knellwolf et al. (2011) and Perić, et. al. (2017). Cossel (2019) provided the solutions for axial displacement, strain and stress in a semi-floating energy pile embedded in a single soil layer that is underlain by a bedrock or another soil layer. The pile tip is flush with the bottom of the relevant soil layer. The solutions for displacement ($u$), strain ($\varepsilon$) and stress ($\sigma$) in a semi-floating energy pile subjected to a thermo-mechanical load in the absence of head restrain are given by Eqs. 1, 2 and 3.

$$u(x) = \frac{\alpha \Delta T \sinh[\psi(x-x_0)]}{\psi \cosh[\psi(L-x_0)]} + \frac{F[E\psi \cosh(\psi x)+k_b \sinh(\psi x)]}{AE\psi[E\psi \sinh(\psi L)+k_b \cosh(\psi L)]} \qquad (1)$$

$$\varepsilon(x) = \frac{\alpha \Delta T \cosh[\psi(x-x_0)]}{\cosh[\psi(L-x_0)]} + \frac{F[E\psi \sinh(\psi x)+k_b \cosh(\psi x)]}{AE[E\psi \sinh(\psi L)+k_b \cosh(\psi L)]} \qquad (2)$$

$$\sigma(x) = E\alpha \Delta T \left\{ \frac{\cosh[\psi(x-x_0)]}{\cosh[\psi(L-x_0)]} - 1 \right\} + \frac{F[E\psi \sinh(\psi x)+k_b \cosh(\psi x)]}{A[E\psi \sinh(\psi L)+k_b \cosh(\psi L)]} \qquad (3)$$

where positive $x$ coordinate is directed upward and it originates at the pile tip. Cossel (2019) also provided the solution for the location of a thermal null point ($x_0$). It is given by Eq. 4.

$$x_0 = \frac{1}{\psi} \tanh^{-1} \left[ \frac{\cosh(\psi L)-1}{\sinh(\psi L)+\frac{k_b}{E\psi}} \right] \qquad (4)$$

The elastic modulus and the coefficient of thermal expansion of the pile are denoted by $E$ and $\alpha$ respectively, while $\Delta T$ is a temperature difference of the pile relative to the surrounding soil. $\Delta T$ is positive in case of heating, and negative in case of cooling. Mechanical load in the form of axial force $F$ is positive in case of a tensile force and negative in case of a compressive force. The parameter $\psi$ contains the pile geometry and relative stiffness of the soil with respect to the pile. It is given by

$$\psi^2 = \left(\frac{p}{A}\right)\left(\frac{k_s}{E}\right) \tag{5}$$

where $p$ and $A$ are the perimeter and the cross-sectional area of the energy pile respectively. The length of the pile is denoted by $L$. The stiffness of the continuous shear spring attached along the pile shaft is denoted by $k_s$, and the stiffness of the normal spring attached at the pile tip is denoted by $k_b$.

In the special case of a nearly ideally fully floating energy pile whereby the restraint at the pile tip is arbitrary small, the thermal null point is located arbitrarily close to the mid-length of the pile. It is noted that ideally fully floating pile would be infinitely long. Nevertheless, a nearly fully floating pile has a finite length and its response is arbitrarily close to that of an ideally fully floating pile. Specifically, the stiffness of the spring attached to the tip ($k_b$) of a nearly fully floating pile can be determined by selecting a desired ratio of mechanically induced axial stresses at the pile tip and pile head. Therefore, substituting $x_0=L/2$, and $k_b=0$ in Eqs. 1 to 3 the following equations are obtained for a nearly fully floating pile.

$$u(x) = \frac{\alpha \Delta T \sinh[\psi(x-L/2)]}{\psi \cosh\psi(L/2)} + \frac{F\cosh(\psi x)}{AE\psi \sinh(\psi L)} \tag{6}$$

$$\varepsilon(x) = \frac{\alpha \Delta T \cosh[\psi(x-L/2)]}{\cosh\psi(L/2)} + \frac{F\sinh(\psi x)}{AE\sinh(\psi L)} \tag{7}$$

$$\sigma(x) = E\alpha\Delta T\left\{\frac{\cosh[\psi(x-L/2)]}{\cosh\psi(L/2)} - 1\right\} + \frac{F\sinh(\psi x)}{A\sinh(\psi L)} \tag{8}$$

## MATERIAL PROPERTIES

In order to conduct thermo-mechanical analysis, a homogenous soil layer, representing the soil layer denoted by A1 by Laloui et al. (2006) is selected. This is one of the four layers that are present within the layered soil profile below Swiss Federal Institute of Technology in Lausanne, Switzerland. In the present study, the soil layer A1 surrounds the energy pile along its entire depth. The parameters of the analytical model used in this study are listed in Table. 1 based on Laloui et al. 2006 and Knellwolf et al. (2011).

Table 1. Parameters of the analytical model

| Parameters | Values |
|---|---|
| Shear spring stiffness, $k_s$ (MPa/m) | 16.7 |
| Energy pile length, $L$ (m) | 26 |
| Energy pile diameter, $D$ (m) | 1 |
| Coefficient of thermal expansion, $\alpha$ (1/°C) | $1 \times 10^{-5}$ |
| Elastic Modulus, $E$ (GPa) | 29.2 |

**PREDICTIONS OF THE ANALYTICAL MODEL**

Four different loading scenarios are considered herein. They can be described as
- Case (*i*): $F < 0$ (compression) and $\Delta T < 0$ (cooling)
- Case (*ii*): $F < 0$ (compression) and $\Delta T > 0$ (heating)
- Case (*iii*): $F > 0$ (tension) and $\Delta T < 0$ (cooling)
- Case (*iv*): $F > 0$ (tension) and $\Delta T > 0$ (heating)

In case of loading scenarios (*ii*) and (*iii*), in which $\Delta T$ and $F$ have different signs, thermal and mechanical loads induce stresses having the same sign. Specifically, in the case of loading scenario (*ii*) both, mechanical and thermal stresses compressive while for loading scenario (*iii*) both, mechanical and thermal are tensile. On the contrary, in case of loading scenarios (*i*) and (*iv*) thermally and mechanically induced stresses have different signs.

First, load scenarios (*i*) and (*ii*) corresponding to the compressive load are considered. In accordance with Knellwolf et al. (2011) an axial compressive force having magnitude 1 MN is applied to the pile head, thus representing a load from the four-storey building. Furthermore, to represent a two-storey building a compressive force having magnitude of 0.5 MN is also applied.

Perić et al. (2020) defined the equivalent thermal load $\Delta T_{eq}$, thus enabling a direct conversion between thermal and mechanical loads. It is given by

$$\Delta T_{eq} = \pm \frac{F}{AE\alpha} \qquad (9)$$

where $A$ is the cross sectional area of the pile, and plus sign applies to load scenarios (*i*) and (*iv*) while minus sign applies to the scenarios (*ii*) and (*iii*). For a mechanical load of magnitude $|F|$, the corresponding magnitude of equivalent thermal load is $|\Delta T_{eq}|$. Eq. 9 is based on the fact that both, mechanical and equivalent thermal load produce displacement and strain of equal magnitude throughout the entire length of an end bearing pile. Thus, for the mechanical load $F = -1$ MN and $F = -0.5$ MN Eq. 9 results in the magnitude of the equivalent thermal loads being equal to 4.36 °C and 2.18 °C respectively. Consequently, for the selected thermal load of ±10 °C the magnitude of actual thermal load is 2.29 times larger than that of the equivalent thermal load whereby the actual mechanical load is -1 MN. In the case that mechanical load is -0.5 MN the magnitude of the equivalent thermal load is 4.58 times larger than that of the equivalent thermal load.

Figures 1 and 2 show thermal, mechanical and combined responses of the fully floating energy pile subjected to load scenario (*i*) and compressive axial forces of 0.5 MN and 1 MN respectively. Figures. 1-a and 2-a show that the thermal null point in both cases is located at the mid-length of the pile as a result of absence of head ant tip restraints in the fully floating energy pile as well as of constant shear spring stiffness throughout the pile length. Nevertheless, the application of the additional mechanical compressive load shifts the initial null point, thus resulting in the combined null point being located lower than the thermal null point. This downward relocation is due to the opposite directions of thermal and mechanical displacements in the bottom part of the pile.

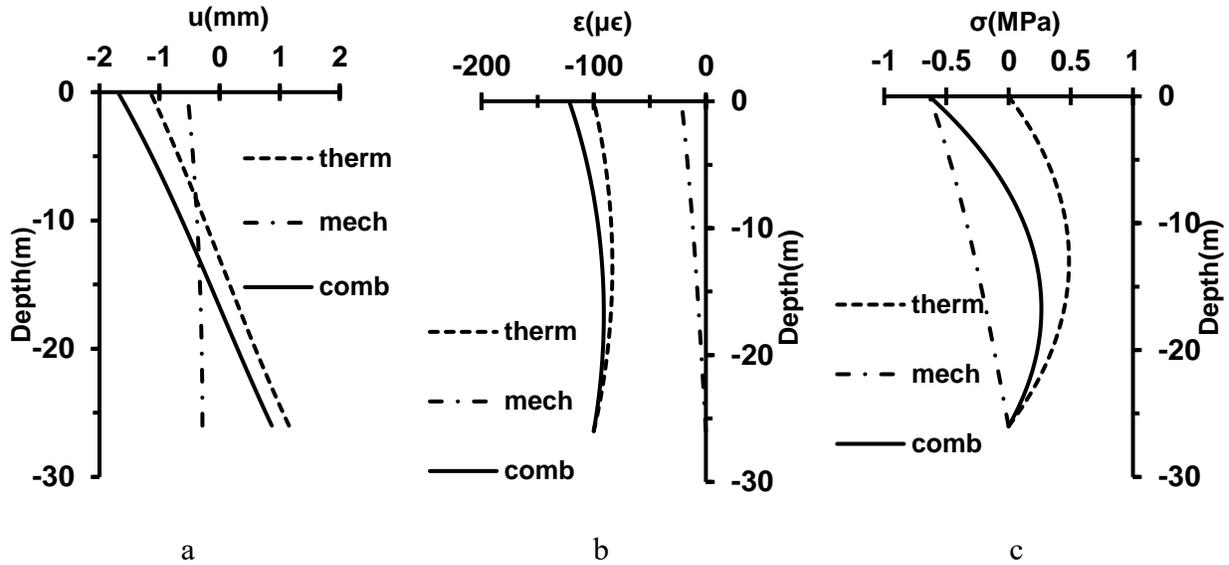

Figure 1. a) Displacement, (b) strain, and (c) stress in a fully floating energy pile subjected to loading scenario (*i*) (*F* = -0.5 MN, and *ΔT*=-10 ºC)

Nevertheless, both, thermal and mechanical loads, induce compressive strains. Furthermore, thermally induced displacement, strain and stress are dominant for the most part because magnitude of the thermal load is significantly larger than that of the mechanical load. It is noted that cooling generates tensile stress in the energy pile while compressive mechanical load induces compressive stress. It is because cooling is dominant as compared to the compressive force that the combined stress is tensile along certain length of the pile (Figures 1c and 2c). The smaller the compressive load the larger the tensile stress and tension zone. Furthermore, the smaller the compressive load the smaller change in the location of a combined null point as compared to the thermal null point.

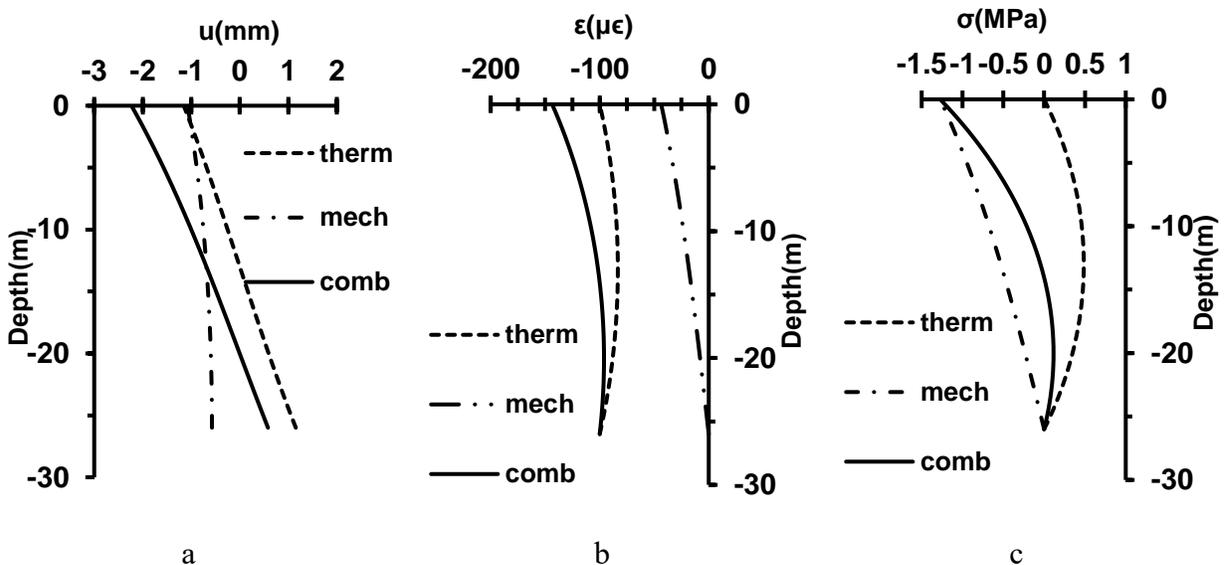

Figure 2. a) Displacement, (b) strain, and (c) stress in a fully floating energy pile subjected to loading scenario (*i*) (*F* = -1MN and *ΔT*= -10 ºC)

Responses of the fully floating energy pile subjected to load scenario (*ii*) whereby the magnitudes of compressive mechanical load are 0.5 MN and 1 MN respectively are depicted in Figures 3 and 4. Heating induces upward displacement in the upper part of the energy pile and downward displacement in its lower part, the thermal null point is located at the mid-depth of the pile. The combined null point in case of heating and compression is located in the upper half of the pile, and it moves further up for the larger magnitude of mechanical load. Heating induces tensile strain, while mechanical strain is compressive as depicted in Figures 3b and 4b. Based on displacement and strain responses thermal load is clearly dominant in the combined response. In the case of load scenario (*ii*) thermally and mechanically induced stresses are both compressive, and they add up in case of a combined load. Thermally induced stress is dominant in the combined response in the lower portion of the pile because thermal stress exceeds the mechanical stress in this region.

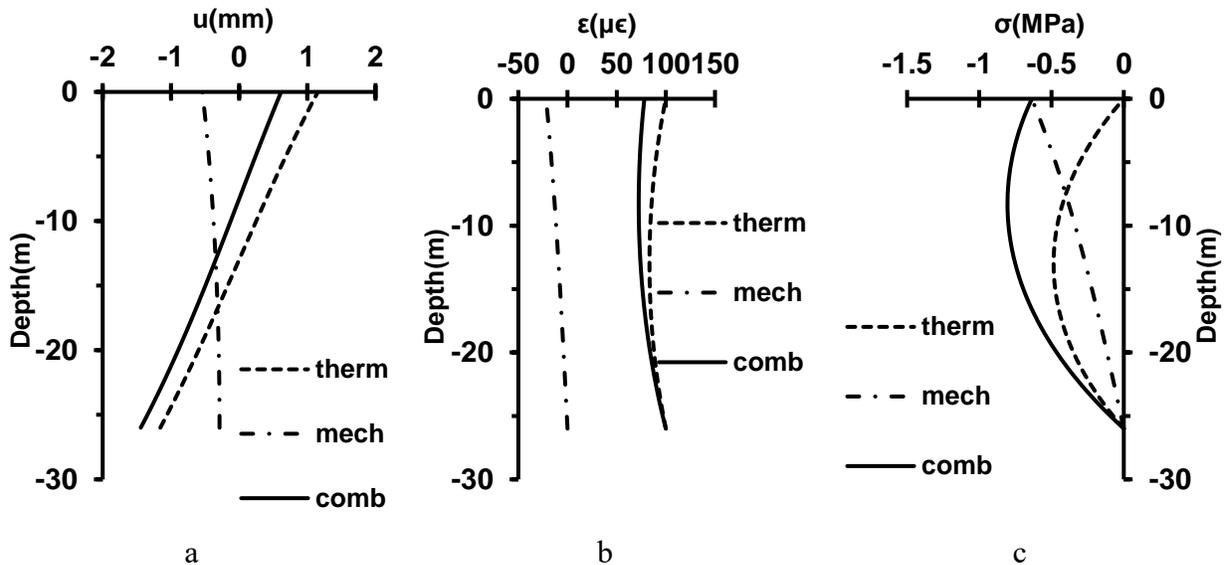

**Figure 3. a) Displacement, (b) strain, and (c) stress in a fully floating energy pile subjected to loading scenario (*ii*) (F = -0.5 MN and *ΔT*=+10 ºC)**

Using Eq. 10 the actual thermal load (*ΔT*) can be expressed in terms of the equivalent thermal load as follows

$$\Delta T = \eta \Delta T_{eq} \tag{10}$$

It is noted that the coefficient $\eta$ is always positive. The larger the values of $\eta$ the more dominant the thermal load is compared to the mechanical load. Figure 5a depicts changes in the location of a combined null point due to increase in $\eta$ for loading scenario (*i*). For $\eta=1$ the combined displacement is directed downward throughout the entire length of the pile, and a combined null point is non-existent. This results in a lack of tensile stress. With increase in $\eta$ thermal response gradually becomes increasingly more dominant. This results in emergence of the combined null point near the pile tip and its gradual movement towards the mid-length of the pile.

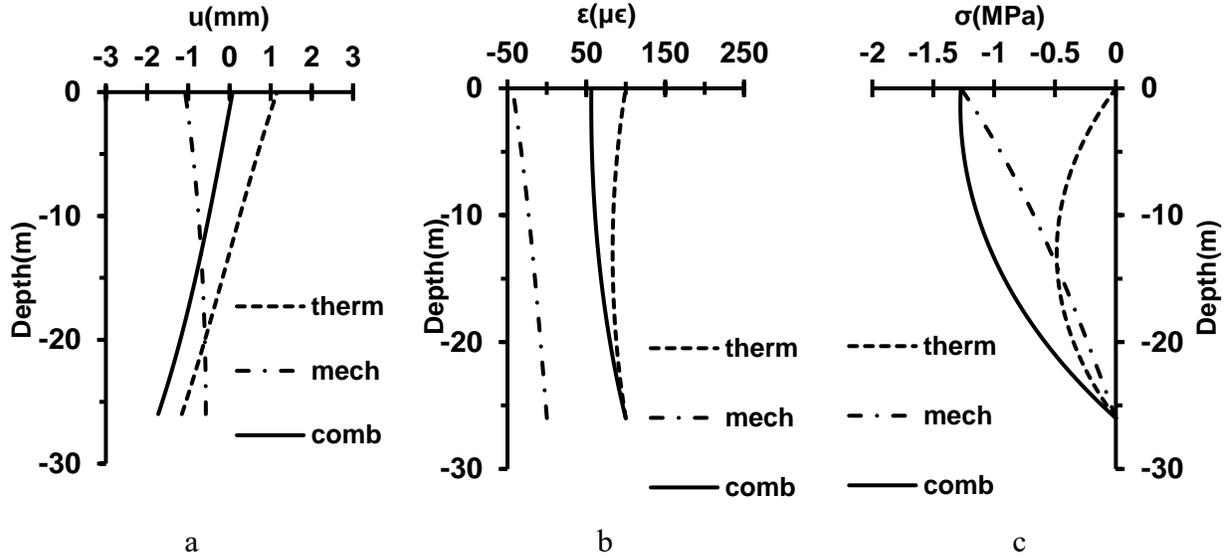

**Figure 4. a) Displacement, (b) strain, and (c) stress in a fully floating energy pile subjected to loading scenario (*ii*) (*F* = -1 MN and *ΔT*=+10 C)**

At the same time, tension zone emerges and it extends upward from the pile tip. As the tension zone becomes longer the magnitude of tensile stress increases. Figures 5 and 6 point out that the location of a combined null point corresponds to the location of the maximum magnitude of the combined stress, which in turn requires perhaps the most important consideration in design of energy piles.

Figure 6a demonstrates development of the combined null point for load scenario (*ii*) and different *η* values. In this case the entire pile moves upward, and the combined null point does not exist for *η* < 2.14. As magnitude of thermal load increases compared to the magnitude of mechanical load the combined null point starts to emerge and gradually move downward towards the mid-length of the pile. As expected, the location of a combined null point coincides with the location of the maximum magnitude of the combined stress (Figure 6b). A larger *η* value indicates a larger dominance of thermal response in the combined response and thus the larger magnitude of a combined compressive stress.

The location of a combined null point ($\bar{x}_0$) can be obtained by combining Eqs. 6, 9 and 10. The solution for load scenarios (*i*) and (*iv*) is given by

$$\bar{x}_0 = \frac{1}{\psi}\tanh^{-1}\left[tanh\left(\frac{\psi L}{2}\right) - \frac{1}{\eta sinh(\psi L)}\right] \tag{11}$$

with limitation

$$\eta > \frac{1}{sinh(\psi L)(tanh(\frac{\psi L}{2}))} \tag{12}$$

The solution for loading scenarios (*ii*) and (*iii*) is given by

$$\bar{x}_0 = \frac{1}{\psi}\tanh^{-1}\left[tanh\left(\frac{\psi L}{2}\right) + \frac{1}{\eta sinh(\psi L)}\right] \tag{13}$$

with limitation

$$\eta > \frac{1}{sinh(\psi L)(1-\tanh(\frac{\psi L}{2}))} \tag{14}$$

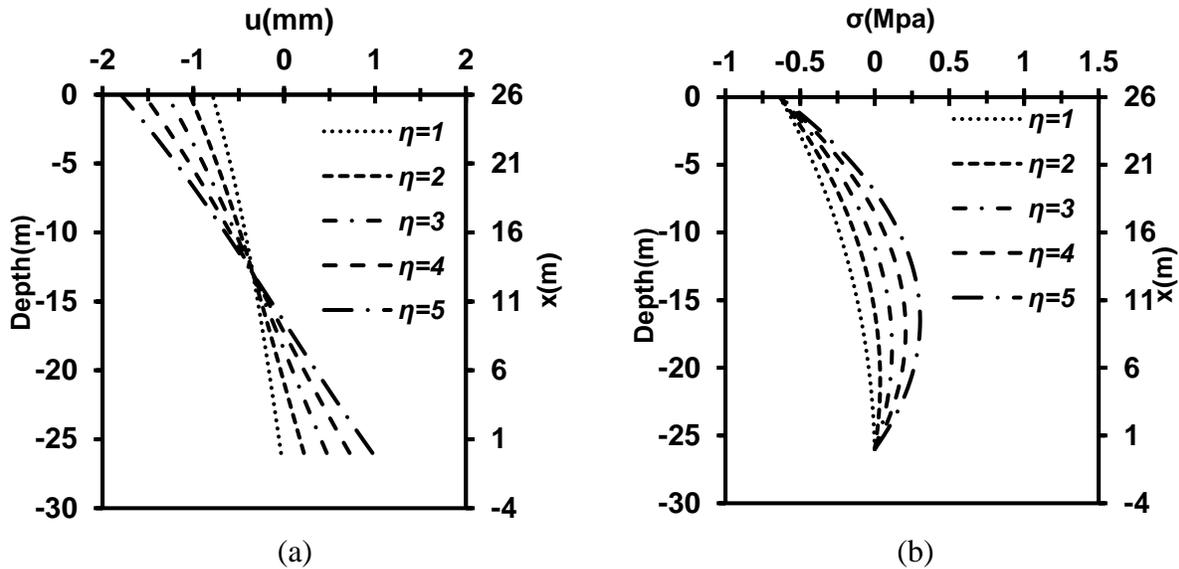

(a)                                 (b)

**Figure 5. a) Displacement, and (b) stress, in a fully floating energy pile subjected to loading scenario (*i*) (*F* = -0.5 MN) and different magnitude of thermal loads**

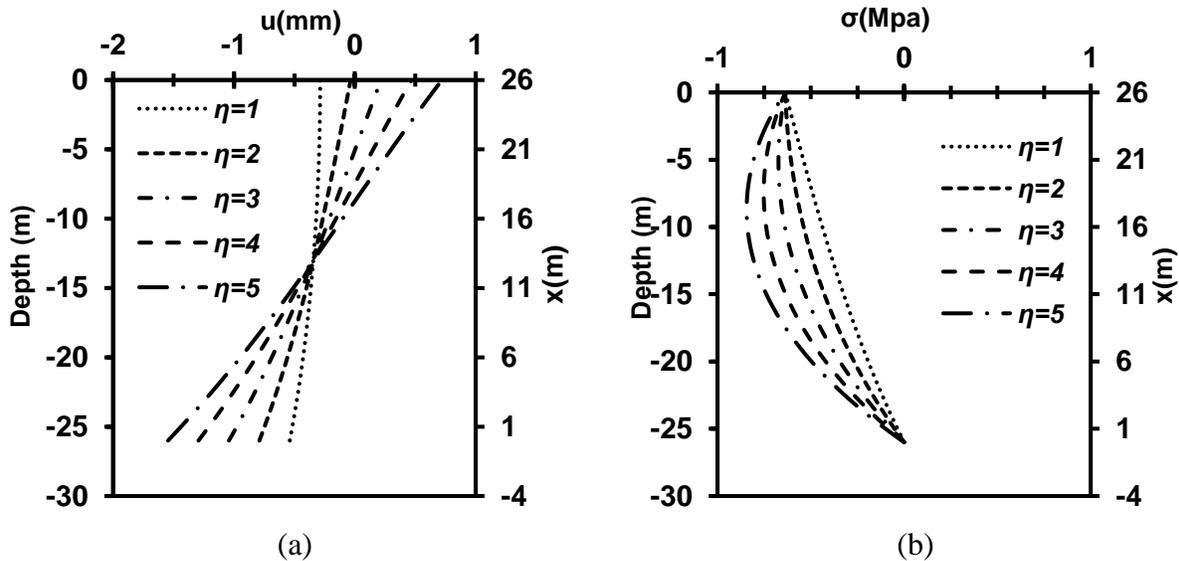

(a)                                 (b)

**Figure 6. a) displacement, and (b) stress, in a floating energy pile subjected to loading scenario (*ii*) (F = -0.5 MN) and different magnitude of thermal loads**

Emergence and progression of the combined null point based on Eqs. 11 and 13 is depicted in Figure 7. For load scenarios (*i*) and (*iv*), where thermal and mechanical loads have the same sign, the combined null point emerges at the pile tip and gradually moves upward towards the mid-depth of the pile with increasing magnitude of thermal load. Whereas, for load scenarios (*ii*) and (*iii*) whereby thermal and mechanical loads have opposite signs, the combined null point emerges at the pile head and gradually moves toward the mid-length of the pile with increasing magnitude of thermal load. In both cases the rate of change of the location of a combined null point is faster at lower values of $\eta$. It is also noted that the curve depicting load scenarios (*i*) and (*iv*) is not the mirror image of the curve depicting load scenarios (*ii*) and (*iii*). Specifically, for load scenarios (*i*) and (*iv*) the combined null point emerges at a smaller thermal load than for load scenarios (*i*) and (*iii*).

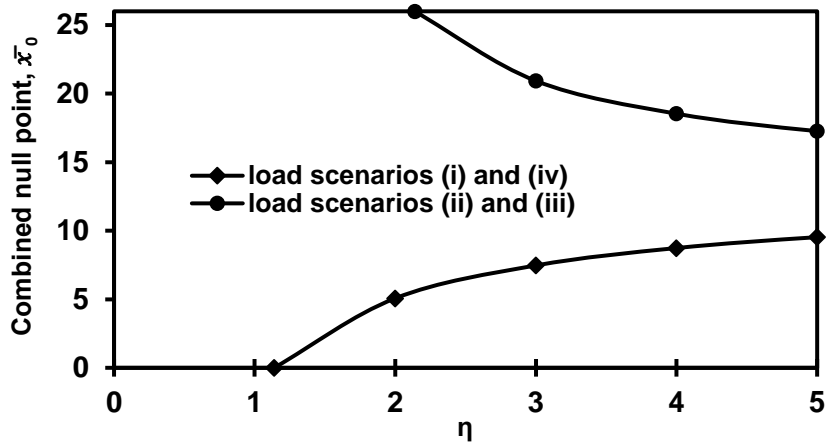

**Figure 7. Location of a combined null point for fully floating energy pile and different magnitudes of thermal load**

**CONCLUSION**

Analytical solutions were implemented to evaluate the response of a fully floating energy pile subjected to different combinations of thermal and mechanical loads. It was found that location of the combined thermal point coincides with the location of the maximum magnitude of thermo-mechanical stress, which is vitally important for design. Mechanical load is typically approximately constant while thermal load changes significantly during the operation of energy piles. Thus, the range of locations of maximum stress magnitude that depends on the operational temperature range can be easily obtained. Furthermore, it is shown that the thermal null point is always located arbitrarily close to the mid-depth of the nearly ideally fully floating energy pile. Nevertheless, the location of a combined null point depends on the magnitude and sign of thermal load as compared to the mechanical load. In case that magnitude of thermal load is small the combined null point is absent. With increasing magnitude of thermal load as compared to mechanical load the combined thermal point starts to emerge and it gradually moves towards the mid-length of the pile. For load combinations where the signs of thermal and mechanical loads are different the combined null point emerges at the pile head and moves downward. For load combinations where the signs of thermal and mechanical loads are the same the combined null point emerges at the pile tip and moves upward. Furthermore, in the latter case the tension zone forms along the portion of the pile length for sufficiently large magnitudes of thermal load.